\documentclass[12pt,a4paper]{article}
\usepackage{amsmath}
\usepackage{amssymb}
\usepackage{amsfonts}
\usepackage{amsthm}
\usepackage{amsmath}
\usepackage{graphicx}
\newtheorem{theorem}{Theorem}

\newtheorem{lemma}{Lemma}
\newtheorem{definition}{Definition}
\newtheorem{remark}{Remark}
\begin{document}


\begin{center}
\Large \bf Homogenization of the Signorini boundary-value\\ problem
in a thick plane junction
\end{center}

\begin{center}
\large \sf  Yulija A.~Kazmerchuk and Taras A.~Mel'nyk
\end{center}

\begin{center}
\small Department of Mathematical Physics, National Taras Shevchenko University of Kyiv,\\
 01033 Kyiv, Ukraine\\
E-mail: ulykaz@ukr.net, \ melnyk@imath.kiev.ua
\end{center}

\begin{abstract} We consider a mixed boundary-value problem for the Poisson equation
in a plane thick junction $\Omega_{\varepsilon}$ which is the union of
a domain $\Omega_0$ and a large number of $\varepsilon-$periodically situated thin rods.
The nonuniform Signorini conditions are given on the vertical sides of the thin rods. The asymptotic analysis of this problem is made as $\varepsilon\to0,$ i.e., when the number of the thin
rods infinitely increases and their thickness tends to zero. With the help of the integral identity method
we prove the convergence theorem and show that the nonuniform Signorini conditions
are transformed (as $\varepsilon\to0)$ in the limiting variational inequalities in the region that is filled up by the thin rods in the limit passage. The existence and uniqueness of the solution to this non-standard limit problem is established.  The convergence of the energy integrals is proved as well.
\end{abstract}


\bigskip
\noindent
{\bf MOS subject classification:} \ 35B27, 35R45, 35J20, 74K30

\section{Introduction and statement of the problem}

Boundary-value problems in thick junctions are mathematical models of
widely used engineering and industrial constructions as well as many other
physical and biological systems with very distinct characteristic
scales. In recent years,  a rich collection of new results on asymptotic analysis of boundary-value problems in thick multi-structures is appeared (see \cite{B-Gaudiello-Melnyk}--\cite{M-V-2}).

In this paper we
homogenize the Signorini problem in a thick plane junction of type $2:1:1$
using the integral identity method developed in \cite{UMJ,M-steklov}.

A thick junction (or thick multi-structure) of type $k:p:d$ \ is the union of some
domain in $\mathbb{R}^{n},$ which is called the junction's body, and a large number of $\varepsilon $-periodically situated thin domains along some manifolds on the boundary of
 the junction's body (see Fig.~\ref{fig1}). This manifold is called the joint zone. Here $\varepsilon $ is a small parameter, which
characterizes the distance between neighboring thin domains and their
thickness.  The type $k:p:d$ of a thick junction refers
respectively to the limiting dimensions of the body, the joint zone, and
each of the attached thin domains.

This classification of thick junctions was given in \cite{MN94}--\cite{M-N-AiA} and \cite{UMJ,M-steklov}, where
rigorous mathematical methods were developed (homogenization, approximation, asymptotic expansions) for analyzing the main boundary-value problems in thick junctions of different types. It was pointed
out that qualitative properties of solutions essentially depend on the
junction type and on the conditions given on the boundaries of the attached
thin domains. In addition, as it was shown in \cite{Fleury-SanchezPal} such problems lose the coercitivity as
$\varepsilon \to 0$ and this creates special difficulties in the
asymptotic investigation. It should be noted that papers \cite{Khrus, Khrus-Kot} are the first papers in this direction.

For the first time a problem known now as the Signorini problem was considered by  Signorini himself in \cite{Signorini}. The sense of the Signorini boundary condition consists in a priori ignorance which of the boundary conditions (Dirichlet or Neumann) are satisfied and where. Many interesting problems in applied mathematics involve the Signorini boundary conditions. Applications arise in groundwater hydrology, in plasticity theory, in crack theory,
in optimal control problems, etc. (see \cite{Kinder-Stam}). Such of these problems as can be recast as variational inequalities become relatively easy to study (see \cite{Kinder-Stam}--\cite{Lions}). Asymptotic investigations of variational inequalities in perforated domains were made in \cite{Zhikov-1}--\cite{Sandrakov-2}.

The results of this preprint will be published in \cite{Mel-Kaz}.

\subsection{Statement of the problem}

\begin{figure}[htbp]\label{fig1}
\begin{center}
\includegraphics[width=8cm]{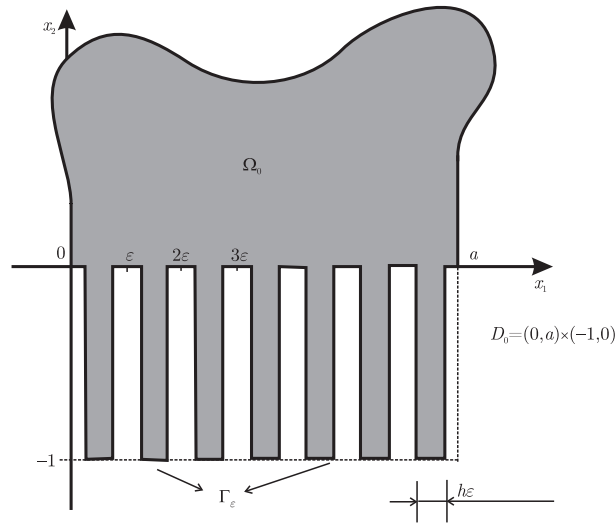} 
\caption{A model plane thick junction $\Omega_{\varepsilon}$ of type 2:1:1}
\end{center}
\end{figure}

Let $a, l$ be positive numbers, h be a fixed number from the interval $(0,1),$ and $N$ be a large positive
integer. Define a small parameter $\varepsilon= \frac{a}{N}.$
A model plane thick junction $\Omega_{\varepsilon}$ (see Fig.~\ref{fig1}) consists of the junction's body
$$
\Omega_{0}=\left\{ x=(x_{1},x_{2}) \in \mathbb{R}^{2}:
\quad  0<x_{1}<a,\quad  0<x_{2}<\gamma(x_{1})\right\},
$$
where $\gamma \in C^{1}([0, a])$, and a large number of the thin rods
$$
G_{j}(\varepsilon)= \left\{ x: \, \left| \frac{x_{1}}{\varepsilon}-
(j+\frac{1}{2})\right|<\frac{h}{2}, \quad x_{2}  \in [-l, 0]
\right\},\quad j=0, 1, ..., N-1,
$$
i. e. $\Omega_{\varepsilon}=\Omega_{0}\cup G_{\varepsilon}$,\quad where
$G_{\varepsilon}=\bigcup \limits_{j=0} \limits^{N-1}G_{j}(\varepsilon).$

The discrete parameter $\varepsilon$ characterizes the distance between the rods and their thickness that
is equal to $\varepsilon h$.  Obviously,  the thin rods fill out the rectangle $D_0 = (0, a)\times (-l, 0)$ in the limit passage as $N\to +\infty \ (\varepsilon\to0).$

Let us denote  the union of vertical sides of the thin rods
$G_{\varepsilon}$ by $S_{\varepsilon};$ \ the union of bases of the thin rods
will be denoted by $\Gamma_{\varepsilon}$.

In $\Omega_{\varepsilon}$ we consider the following boundary value
problem
\begin{equation}\label{1.1}
\left\{
  \begin{array}{c}
     -\Delta u_{\varepsilon}(x)=f(x),\quad x \in \Omega_{\varepsilon},
\\\\
          u_{\varepsilon}(x)\leq g(x),\quad \partial_{\nu}u_{\varepsilon}(x)\leq\varepsilon
d(x), \quad x \in S_{\varepsilon},
\\\\
     \left(u_{\varepsilon}(x)-g(x)\right)\left(\partial_{\nu}u_{\varepsilon}(x)-\varepsilon
        d(x)\right)=0, \quad x \in S_{\varepsilon},
\\\\
u_{\varepsilon}(x)=0, \quad x \in \Gamma_{\varepsilon},\qquad
     \partial_{\nu}u_{\varepsilon}(x)=0,\quad x \in
        \partial\Omega_{\varepsilon}\setminus(S_{\varepsilon}\cup\Gamma_{\varepsilon}),
  \end{array}
\right.
\end{equation}
where $\partial_{\nu}=\partial/\partial \nu$ is the outward normal
derivative, \ $f, g, d$ are given functions.

We assume that  $f \in L^{2}(\Omega_{1}),$ where
$\overline{\Omega_{1}}=\overline{\Omega_{0}}\cup\overline{D_{0}},$
the function $d$ belongs to the Sobolev space $H^{1}({D_{0}}),$ and
$$
g\in H^{1}(D_{0}; \, I_{l}\cup I_{0}) =\left\{ v \in
H^{1}(D_{0}):\quad \left.v\right|_{I_{l}\cup I_{0}}=0\right\},
$$
where
 $I_{l}=\left\{x: \ \ x_{1} \in (0,a),\quad x_{2}=-l\right\}, \quad I_{0}=\left\{x: \ \  x_{1} \in (0,a),\quad
x_{2}=0\right\}.$

\medskip

{\sf Our goal is to study the asymptotic behavior of the solution $u_\varepsilon$ to problem (\ref{1.1}) as $\varepsilon\to0,$ i.e., when the number of the thin rods infinitely increases and their thickness tends to zero.
}

\section{Definitions of the weak solution and its existence}

In the Sobolev space $H^{1}(\Omega_\varepsilon; \, \Gamma_\varepsilon) =\left\{ u \in
H^{1}(\Omega_\varepsilon) : \ \left.u\right|_{\Gamma_\varepsilon}=0\right\},
$ we define subset
$$
K_{\varepsilon}=\{\varphi \in H^{1}(\Omega_{\varepsilon}; \Gamma_{\varepsilon}) : \quad
\left.\varphi\right|_{S_\varepsilon}\leq \left.g \right|_{S_\varepsilon} \ \text{a.e. on} \ S_{\varepsilon}
\},
$$
where $\left.\psi\right|_{S}$ denotes the trace of a Sobolev function $\psi$ on a curve $S.$ Obviously, $K_\varepsilon$ is a closed and convex set for every fixed value of $\varepsilon.$

Let us suppose the existence of a classical solution to problem (\ref{1.1}). We can regard that $g=0$ in $\Omega_{0}.$ Multiplying the equation of problem (\ref{1.1}) by the function $(u_{\varepsilon}-g),$ \
integrating by parts in $\Omega_{\varepsilon}$ and taking into account the boundary conditions for $u_{\varepsilon},$
we obtain
\begin{equation}\label{2.1}
\int\limits_{\Omega_{0}}\nabla u_{\varepsilon} \cdot \nabla
    u_{\varepsilon}\,dx+ \int\limits_{G_{\varepsilon}}\nabla
    u_{\varepsilon} \cdot \nabla (u_{\varepsilon}- g)\,dx=
    \int\limits_{\Omega_{0}} f u_{\varepsilon}\,dx +
    \int\limits_{G_{\varepsilon}} f
    (u_{\varepsilon}-g)\,dx+\varepsilon\int\limits_{ S_{\varepsilon}}
    d(x)(u_{\varepsilon}-g)\, ds,
\end{equation}
where $\nabla v \cdot \nabla w = \sum\limits_{j=1}\limits^{n}
\frac{\partial v}{\partial x_{j}}\frac {\partial w}{\partial
x_{j}}.$

Now we take any function $\varphi\in K_{\varepsilon}$ and multiply the equation of the problem (\ref{1.1}) by
$(\varphi-g).$ Similar as before we get
$$
\int\limits_{\Omega_{0}}\nabla
u_{\varepsilon} \cdot \nabla \varphi \,
    dx+ \int\limits_{G_{\varepsilon}}\nabla u_{\varepsilon} \cdot \nabla
    (\varphi - g)\,dx=
$$
\begin{equation}\label{2.2}
    =\int\limits_{\Omega_{0}} f  \varphi \, dx +
    \int\limits_{G_{\varepsilon}} f (\varphi-g) \,dx+\varepsilon
    \int\limits_{ S_{\varepsilon}} d(x) (\varphi-g)\, ds +\int\limits_{
    S_{\varepsilon}}(\partial_{\nu} u_{\varepsilon}- \varepsilon
    d(x))(\varphi-g) \,ds.
\end{equation}
Since  $\partial_{\nu}u_{\varepsilon}(x)\leq\varepsilon
d(x)$ and $\varphi(x)\leq g(x)$ a.e. in $S_{\varepsilon},$
\begin{equation}\label{2.3}
   \int\limits_{ S_{\varepsilon}}(\partial_{\nu} u_{\varepsilon}-
    \varepsilon d(x))(\varphi-g)\, ds \geq 0.
\end{equation}
Taking into account (\ref{2.3}), it   follows from equality (\ref{2.2}) that
\begin{equation}\label{2.4}
\int\limits_{\Omega_{0}}\nabla u_{\varepsilon}\cdot\nabla \varphi\,
dx+ \int\limits_{G_{\varepsilon}}\nabla u_{\varepsilon}\cdot\nabla
(\varphi- g)\,dx \geq \int\limits_{\Omega_{0}} f \varphi \,dx +
\int\limits_{G_{\varepsilon}} f (\varphi-g)\,dx+\varepsilon
\int\limits_{ S_{\varepsilon}} d(x) (\varphi-g) \,ds
\end{equation}

\begin{definition}\label{def-1}
A function  $u_{\varepsilon} \in K_{\varepsilon}$ is called a weak solution to problem
(\ref{1.1}) if it satisfies the integral equality (\ref{2.1}) and integral
inequality (\ref{2.4}) for arbitrary function $ \varphi \in K_{\varepsilon}.$
\end{definition}

Another definition is as follows.

\begin{definition}\label{def-2}
 A function  $u_{\varepsilon} \in
K_{\varepsilon}$ is called a weak solution to problem (\ref{1.1}) if it satisfies the integral inequality
\begin{equation}\label{ozn2.5}
\int\limits_{\Omega_{\varepsilon}}\nabla u_{\varepsilon} \cdot
\nabla (\varphi - u_{\varepsilon})\,dx \geq
\int\limits_{\Omega_{\varepsilon}} f (\varphi-u_{\varepsilon}) \,dx
+\varepsilon \int\limits_{ S_{\varepsilon}} d(x)
(\varphi-u_{\varepsilon}) \,ds \quad \forall\  \varphi \in K_{\varepsilon}.
\end{equation}
\end{definition}

Let us show that these definitions are equivalent. Subtracting
equality (\ref{2.1}) from inequality (\ref{2.4}), we arrive at (\ref{ozn2.5}).
Setting $\varphi =\left\{
\begin{array}{ll}
  0 , & x\in  \Omega_0 , \\
  g, & x\in  G_\varepsilon ,
\end{array}\right.$ in (\ref{ozn2.5}),
we have
\begin{equation}\label{2.6}
-\int\limits_{\Omega_{0}}|\nabla u_{\varepsilon}|^2 \,dx+ \int\limits_{G_{\varepsilon}}\nabla
    u_{\varepsilon}\cdot\nabla (g- u_{\varepsilon})\,dx
    \geq-\int\limits_{\Omega_{0}} f u_{\varepsilon}\, dx +
    \int\limits_{G_{\varepsilon}} f (g-u_{\varepsilon})\,dx+\varepsilon
    \int\limits_{ S_{\varepsilon}} d(x) (g-u_{\varepsilon}) \,ds.
\end{equation}
Putting $\varphi=\left\{
\begin{array}{ll}
  2u_{\varepsilon} , & x\in  \Omega_0 , \\
  2u_{\varepsilon}-g, & x\in  G_\varepsilon ,
\end{array}\right.$ in (\ref{ozn2.5}), we get the reversed inequality
\begin{equation}\label{2.7}
\int\limits_{\Omega_{0}}|\nabla u_{\varepsilon}|^2 \,dx+ \int\limits_{G_{\varepsilon}}\nabla
    u_{\varepsilon}\cdot\nabla (u_{\varepsilon}-g)\,dx
    \geq\int\limits_{\Omega_{0}} f u_{\varepsilon}\, dx +
    \int\limits_{G_{\varepsilon}} f (u_{\varepsilon}-g)\,dx+\varepsilon
    \int\limits_{ S_{\varepsilon}} d(x) (u_{\varepsilon}-g)\, ds.
\end{equation}
This means that (\ref{2.1}) holds.
 Setting $\varphi=\left\{
\begin{array}{ll}
   \psi+u_{\varepsilon}, & x\in  \Omega_0 , \\
   \psi+u_{\varepsilon}-g, & x\in  G_\varepsilon ,
\end{array}\right.$ in (\ref{ozn2.5}), where $\psi$
is arbitrary function from $K_{\varepsilon},$ we get (\ref{2.4}).

It is well known (see for instance  \cite{Kinder-Stam}--\cite{Lions}) that there exists a unique solution of inequality~(\ref{ozn2.5}) for any fixed value of $\varepsilon.$

\section{ Auxiliary uniform estimates}
To homogenize boundary-value problems in thick multi-structures with
the nonhomogeneous Neumann or Fourier conditions on the boundaries
of the thin attached domains the method of special integral identities was
proposed in  \cite{UMJ,M-steklov}. For our problem this identity is as follows (see \cite[Lemma~1]{M-steklov})
\begin{equation}\label{4.1}
  \frac{\varepsilon h}{2}\int_{S_{\varepsilon} }v
  \: dx_2=\int_{G_{\varepsilon}} v \: dx -
  \varepsilon\int_{G_{\varepsilon} }Y\left(\frac{x_1}{\varepsilon
}\right)\partial_{x_1}v \: dx \qquad \forall\: v\in
H^1\bigl(G_{\varepsilon} \bigr),
\end{equation}
where $Y(\xi )=-\xi +[\xi ]+ \frac12,$ $[\xi]$ is the integral part
of $\xi .$

Using (\ref{4.1}) and taking into account that
$\max_\mathbb{R}|Y|\leq 1,$  we get
\begin{equation}\label{4.2}
\| v\|_{L^2(S_\varepsilon)}\leq C_1 \varepsilon^{-\frac12} \|v
\|_{H^1(G_\varepsilon)}\quad \forall \ v \in H^1(G_\varepsilon).
\end{equation}

\begin{remark}
Here and in  what follows all constants $\{C_{i}\}$ and $\{c_{i}\}$ in
inequalities are independent of the parameter $\varepsilon.$
\end{remark}

Also it will be very important for our research the following uniform estimates.

\begin{lemma}[\cite{ZAA99}]\label{lemma-1}
The usual norm
$
\|u\|_{H^{1}(\Omega_{\varepsilon})}=\left(\int_{\Omega_{\varepsilon}}\left(|\nabla
u|^{2}+u^{2}\right)dx\right)^{\frac{1}{2}}
$
in $H^{1}(\Omega_{\varepsilon}, \Gamma_\varepsilon)$ and a norm
 $\|\cdot\|_{\varepsilon},$ which is generated by scalar product
 $$
(u,v)_{\varepsilon}=\int_{\Omega_{\varepsilon}}\nabla
u\cdot\nabla v \: dx, \quad u,v \in H^{1}(\Omega_{\varepsilon}, \Gamma_\varepsilon),
$$
are uniformly equivalent, i.e., there exist constants $C_{1}>0$
and $\varepsilon _{0}>0$ such for all $\varepsilon \in
(0,\varepsilon _{0})$ and for all $u \in H^{1}(\Omega_{\varepsilon}, \Gamma_\varepsilon)$ the following estimates
\begin{equation}\label{lema-est-4.3}
\|u\|_{\varepsilon } \leq \|u\|_{H^{1}(\Omega_{\varepsilon })} \leq
C_1\|u\|_{\varepsilon}
\end{equation}
hold.
\end{lemma}

\begin{remark}
In fact, in Lemma~\ref{lemma-1} the following Friedrich inequality
\begin{equation}\label{4.4}
\|u\|_{L^2(\Omega_\varepsilon)} \leq C_2 \|\nabla u\|_{L^2(\Omega_\varepsilon)} \quad \forall \ u \in H^{1}(\Omega_{\varepsilon}, \Gamma_\varepsilon)
\end{equation}
was proved.
\end{remark}

Using the Cauchy--Bunyakovsky integral inequality and Cauchy's inequality with $\delta> 0$ $(2ab\leq\delta a^{2}+\delta^{-1}b^{2}$ for any positive numbers $a$ and $b),$ with the help of (\ref{4.2}) and
(\ref{4.4}) we deduce from (\ref{2.1}) that
\begin{multline}\label{4.5}
\int_{\Omega_{\varepsilon}}|\nabla
u_{\varepsilon}|^{2}\,dx\leq  c_{0}(\delta_{1}+ \delta_{2} +\delta_{3})
\, \|\nabla u_{\varepsilon}\|^{2}_{L^{2}(\Omega_{\varepsilon})} +
\\
+ c_1 (1 + \delta_{1}^{-1}) \, \|g\|^{2}_{H^{1}(D_{0})} +
c_2(1+ \delta_{2}^{-1}) \, \|f\|^{2}_{L^{2}(\Omega_{\varepsilon})} +
c_3(1+ \delta_{3}^{-1}) \, \|d\|^{2}_{H^{1}(D_{0})}.
\end{multline}

Choosing  $\delta_{1}, \delta_{2}, \delta_{3}$ so that
$c_{1}(\delta_{1}+\delta_{2}+\delta_{3})<\frac{1}{2},$  we have
\begin{equation}\label{4.6}
\int_{\Omega_{\varepsilon}}|\nabla u_{\varepsilon}|^{2}dx \leq
c_{4}\left(\|f\|^{2}_{L^{2}(\Omega_1)}+\|g\|^{2}_{H^{1}(D_{0})}+\|d\|^{2}_{H^{1}(D_{0})}\right).
\end{equation}

By virtue of (\ref{lema-est-4.3})  we obtain from (\ref{4.6}) the following uniform estimate
\begin{equation}\label{4.7}
\|u_{\varepsilon}\|_{H^{1}(\Omega_{\varepsilon})}\leq C_3.
\end{equation}

\section{Convergence theorem}

In the sequel, $\widetilde{u}$ denotes the zero extension of a function $u$ defined on $G_\varepsilon$
 to the rectangle $D_0,$ which is filled up by the thin rods in the limit passage as $\varepsilon\to0,$
 namely:
\begin{equation*}
\widetilde{u}(x) =\left\{
\begin{array}{ll}
  u(x), & x\in  G_\varepsilon , \\
  0, & x\in D_0\setminus G_\varepsilon.
\end{array}\right.
\end{equation*}

If $u\in H^1(G_\varepsilon, \Gamma_\varepsilon),$ then due to the rectilinearity of the boundaries of the thin rods $\widetilde{u} $ belongs to the anisotropic Sobolev space
\begin{equation}\label{5.1}
H^{0,1}(D_0; I_{l})=\left\{v\in L^2(D_0):\: \exists \ \mbox{
weak derivative } \ \partial_{x_2}v\in L^2(D_0) \ \ \mbox{and} \ \ \left.v\right|_{I_l}=0 \right\}
\end{equation}
and  $\partial_{x_2}\widetilde{u} = \widetilde{\partial_{x_2}u}$ \, a.e. in $D_0.$

\begin{theorem}\label{Th1}
The solution $u_\varepsilon$ to problem (\ref{1.1}) satisfies the relations
\begin{equation}\label{limits}
\left.
\begin{array}{rcll}
\left.u_\varepsilon\right|_{\Omega_0}
&\stackrel{w}{\longrightarrow}&u^+_0 &\mbox{weakly in}\ \ H^1(\Omega _0),
\\[2mm]
\widetilde{u}_\varepsilon &\stackrel{w}{\longrightarrow}& h \, u_0^{-}
    &\mbox{weakly in}\ \ W^{0,1}(D_0, I_{l}),
\\[2mm]
\widetilde{\partial_{x_1}u_\varepsilon}
&\stackrel{w}{\longrightarrow}&0 & \mbox{weakly in}\ \ L^2(D_0)
\end{array} \right\}
\quad \text{as} \quad \varepsilon \to 0,
\end{equation}
and the function $u_0(x) =\left\{
\begin{array}{ll}
  u_0^+ , & x\in  \Omega_0 , \\
  u_0^-,  & x\in D_0 ,
\end{array}\right.$ is a unique solution of the following problem
\begin{equation}\label{limit-problem} \left\{
\begin{array}{cccc}
-\Delta u_0^+(x)&=f(x), & \quad x\in \Omega_0,
\\[2mm]
-h\, \partial^2_{x_2 x_2} u_0^{-}(x)&\leq h\,f(x) + 2\, d(x), & \quad x \in D_0,
\\[2mm]
u_0^-(x)&\leq g(x),  & \quad x \in D_0,
\\[2mm]
(u_0^-(x)-g(x))\, (h\,\partial^2_{x_2 x_2 } u_0^{-}(x)&+hf(x)+2d(x))=0, & \quad x\in D_0,
\\[2mm]
\partial_{\nu}u_0^+(x)&=0,  &\quad x\in \partial\Omega_0\setminus I_0,
\\[2mm]
u_0^{-}(x_1,-l)&=0,& \quad x_1\in [0,a],
\\[2mm]
u_0^+(x_1,0)&=u_0^-(x_1,0), & \quad x_1\in [0,a],
\\[2mm]
\partial_{x_2}u_0^+(x_1,0)& = h\, \partial_{x_2}u_0^-(x_1,0),& \quad x_1\in [0,a],
\end{array}\right.
\end{equation}
which is called the homogenized problem for (\ref{1.1}).

Furthermore, the following energy convergence holds
$$
\lim_{\varepsilon \rightarrow 0} E_{\varepsilon}(u_{\varepsilon})= E_{0}(u_0),
$$
where
$$
E_{\varepsilon}(u_{\varepsilon}) = \int\limits_{\Omega_{0}}|\nabla u_{\varepsilon}|^2\,dx, \quad
E_{0}(u_0)=\int\limits_{\Omega_{0}}|\nabla u_0^+ |^2\, dx + h \int\limits_{D_0} |\partial_{x_2}u_0^-|^2\, dx.
$$
\end{theorem}

Before the proof of Theorem~\ref{Th1} we investigate the homogenized problem
(\ref{limit-problem}).

\subsection{The solvability of the homogenized problem}

We see that the homogenized problem (\ref{limit-problem}) is a non-standard boundary-value problem that consists of the Poisson equation in the junction body $\Omega_0,$  the variational inequalities in $D_0$
and the transmission conditions in the joint zone $I_0.$ Therefore, at first we give the definition of the weak solution to this problem and then with the help of general approach in the theory of variational inequalities we prove the existence and uniqueness.

Let us introduce partly anisotropic Sobolev space
\[
{\cal H}(\Omega_1, D_0; I_l) = \left\{u\in L^2(\Omega_1):
\  \ \exists\, \partial_{x_2} u\in L^2(\Omega_1),\quad  \left.u\right|_{\Omega_0}\in H^1(\Omega_0), \quad
\left.u\right|_{D_0} \in H^{0,1}(D_0; I_l)\right\},
\]
where $H^{0,1}(D_0; I_l)$ is defined in (\ref{5.1}). It follows from properties of anisotropic Sobolev spaces (see \cite{U}) that the traces of the restrictions $u^+:= \left.u\right|_{\Omega_0}$ and $u^-:=\left.u\right|_{D_0}$ on $I_0$ are equal. In addition, since traces of functions from ${\cal H}(\Omega_1, D_0; I_l)$ vanish on $I_l,$ there exists a constant
$C_0$ such that
\[
\int_{\Omega_1} u^2\, dx \le C_0 \Bigl(\int_{\Omega_{0}}|\nabla u^+|^2 \,dx + \int_{D_0} |\partial_{x_2}u^-|^2
\, dx\Bigr)  \quad \forall \ u \in {\cal H}(\Omega_1, D_0; I_l).
\]
In ${\cal H}(\Omega_1, D_0; I_l)$ we introduce a norm $\|\cdot\|_{\cal H}$, which is generated by a
scalar product
$$
(u, v)_{\cal H}=\int_{\Omega_{0}}\nabla u^+ \cdot \nabla v^+ \,dx + h \int_{D_0} \partial_{x_2}u^- \,\partial_{x_2} v^-
\, dx , \quad u,\ v \in {\cal H}(\Omega_1, D_0; I_l).
$$

We now define subset
$
K_0 = \{\varphi \in {\cal H}(\Omega_1, D_0; I_l) : \ \varphi^-\leq g \ \ \text{a.e. in} \ D_0\}
$
in ${\cal H}(\Omega_1, D_0; I_l).$
Obviously, $K_0$ is a closed and convex set.

\begin{definition}\label{def-3}
A function  $u_0 \in K_0$ is called a weak solution of problem (\ref{limit-problem}) if
it satisfies the integral equality
\begin{equation}\label{def3.1}
\int\limits_{\Omega_{0}}\nabla u_0^+ \cdot \nabla u_0^+ \, dx+h
\int\limits_{D_0} \partial_{x_2}u_0^-
\partial_{x_2}(u_0^- - g)\, dx=\int\limits_{\Omega_0}f u_0^+ \, dx +h \int\limits_{D_0} f \,(u_0^-
-g)\, dx+2 \int\limits_{D_0} d \,(u_0^- -g)\, dx ,
\end{equation}
and the integral inequality
\begin{equation}\label{def3.2}
\int\limits_{\Omega_{0}}\nabla u_0^+ \cdot \nabla \varphi \, dx+h
\int\limits_{D_0} \partial_{x_2}u_0^-
\partial_{x_2}(\varphi - g)\, dx\geq\int\limits_{\Omega_0}f \varphi \, dx +h \int\limits_{D_0} f \,(\varphi
-g)\, dx+2 \int\limits_{D_0} d \,(\varphi -g)\, dx
\end{equation}
for arbitrary function   $ \varphi \in K_0.$
\end{definition}

If there exists a classical solution to the homogenized problem (\ref{limit-problem}), then relations (\ref{def3.1}) and (\ref{def3.2}) can be obtained by the same way as relations (\ref{2.1}) and (\ref{2.4}) in Definition~\ref{def-1}.
Similarly as we proved the equivalence of Definitions~\ref{def-1} and ~\ref{def-2}, we can show the equivalence
of Definition~\ref{def-3} to the following definition.

\begin{definition}\label{def-4}
A function  $u_0 \in K_0$ is called a weak solution to problem (\ref{limit-problem}) if
it satisfies the integral inequality
$$
\int_{\Omega_{0}}\nabla u_0^+ \cdot \nabla (\varphi-u_0^+)
\, dx + h \int_{D_0} \partial_{x_2}u_0^-\,
\partial_{x_2}(\varphi - u_0^-)\, dx\geq
$$
\begin{equation}\label{def4.1}
\geq\int_{\Omega_0}f (\varphi-u_0^+) \, dx +h
\int_{D_0} f \,(\varphi -u_0^-)\, dx + 2 \int_{D_0} d \,(\varphi -u_0^-)\, dx
\end{equation}
for arbitrary function   $ \varphi \in K_0$.
\end{definition}

We now give the third definition of the weak solution to problem (\ref{limit-problem}).

\begin{definition}\label{def-5} A function  $u_0 \in K_0$ is called a
weak solution to problem (\ref{limit-problem}) if it satisfies the
integral inequality
$$
\int_{\Omega_{0}}\nabla \varphi \cdot \nabla (\varphi-u_0^+)
\, dx + h \int_{D_0} \partial_{x_2}\varphi\, \partial_{x_2}(\varphi - u_0^-)\, dx\geq
$$
\begin{equation}\label{def5.1}
\geq\int_{\Omega_0}f (\varphi-u_0^+) \, dx + h \int_{D_0} f \,(\varphi -u_0^-)\, dx + 2
\int_{D_0} d \,(\varphi -u_0^-)\, dx
\end{equation}
for arbitrary function   $ \varphi \in K_0$.
\end{definition}

Let us prove that Definition~\ref{def-4} and Definition~\ref{def-5} are equivalent.
Adding  the inequality
$$
\int_{\Omega_0}\nabla (\varphi-u_0^+) \cdot \nabla (\varphi-u_0^+)\, dx + h\int_{D_0}
\partial_{x_2}(\varphi-u_0^-) \, \partial_{x_2}(\varphi-u_0^-)\geq0 \quad (\varphi \in K_0)
$$
to inequality (\ref{def4.1}), we get (\ref{def5.1}).
Now we take any  $\psi \in K_0.$ Setting $\varphi=u_0+t(\psi-u_0) \in K_0$ \ (for any $t\in [0,1])$ in inequality (\ref{def5.1}), we obtain
$$
\int_{\Omega_{0}}\nabla (u_0^+ +t(\psi-u_0^+)) \cdot \nabla (\psi-u_0^+) \, dx + h \int_{D_0} \partial_{x_2}(u_0^-+t(\psi-u_0^-))\, \partial_{x_2}(\psi-u_0^-) \, dx\geq
$$
\begin{equation}\label{5.2}
\geq\int_{\Omega_0}f (\psi-u_0^+) \, dx +h \int_{D_0}
f \,(\psi-u_0^-)\, dx + 2 \int_{D_0} d \,(\psi-u_0^-)\, dx.
\end{equation}
Passing to the limit in (\ref{5.2}) as $t\rightarrow0,$ we arrive at (\ref{def4.1}).
Thus, all Definitions~\ref{def-3}, \ref{def-4} and \ref{def-5} are equivalent.

We can re-write inequality (\ref{def4.1}) in the following form
\begin{equation}\label{5.3}
(u,\varphi- u)_{\cal H} \geq \langle F, \varphi- u\rangle \quad \forall \ \varphi \in K_0,
\end{equation}
where $F$ is a linear continuous functional on ${\cal H}(\Omega_1, D_0; I_l)$ defined by the formulae
\[
\langle F, w \rangle = \int_{\Omega_0} f \,w^+\, dx + h\int_{D_0} f\, w^-\,dx + 2\int_{D_0} d\, w^-\, dx
\quad \text{for all} \ w \in {\cal H}(\Omega_1, D_0; I_l).
\]
Using the theory of variational inequalities in Hilbert spaces (see \cite[Sec.~2]{Kinder-Stam}), we can state that
there exists a unique solution of the inequality (\ref{5.3}) and consequently of the homogenized problem~(\ref{limit-problem}).

\subsection{The proof of Theorem~\ref{Th1}}

\paragraph{1.}
 From (\ref{4.7}) it follows that
$\|u_{\varepsilon}\|_{H^{1}(\Omega_{0})}\leq C_3,$
$\|\widetilde {u_{\varepsilon}}\|_{L^{2}(D_0)}\leq C_3$ and
$\|\widetilde{\partial_{x_i}u_\varepsilon}\|_{L^{2}(D_{0})}\leq
C_3,$ $ i=1,2.$ Therefore we can choose a subsequence $\{\varepsilon
'\}\subset\{\varepsilon \}$ (again denoted by $\varepsilon $),
such that
\begin{equation}\label{5.4}
\left.
\begin{array}{rcll}
\left.u_\varepsilon\right|_{\Omega_0}
&\stackrel{w}{\longrightarrow}&u^+_0 &\mbox{weakly in}\ \ H^1(\Omega _0),
\\[2mm]
\widetilde{u}_\varepsilon &\stackrel{w}{\longrightarrow}& h\,u_0^{-}
    &\mbox{weakly in}\ \ L^2(D_0),
\\[2mm]
\widetilde{\partial _{x_i}u_\varepsilon } &\stackrel{w}{\longrightarrow}&\gamma _i & \mbox{weakly in}\ \ L^2(D_0), \quad i=1,2,
\end{array} \right\}
\quad \text{as} \quad \varepsilon \to 0,
\end{equation}
where $u^+_0,$ $u^{-}_0,$ $\gamma _1,$ $\gamma _2$ are some
functions which will be determined later.

At first we determine  $\gamma _2.$ Take any function $\psi \in
C_0^\infty (D_0)$ and  perform the following calculations
$$
\int\limits_{D_0}\widetilde{\partial_{x_2}u_\varepsilon } \, \psi \:
dx =\int\limits_{D_0}\partial_{x_2}\widetilde{u_\varepsilon} \, \psi
\: dx = \int\limits_{G_\varepsilon}\partial_{x_2}u_\varepsilon \, \psi
\: dx  = -\int\limits_{G_\varepsilon}u_\varepsilon\,
\partial_{x_2}\psi \: dx= -\int\limits_{D_0}\widetilde{u_\varepsilon} \,
\partial_{x_2}\psi \: dx.
$$
Passing to the limit in this identity, as $\varepsilon \rightarrow
0,$ we obtain
\begin{equation}\label{gamma_2}
 \int_{D_0}\gamma _2  \, \psi \: dx =
-h \int_{D_0}u_0 ^{-}\, \partial_{x_2}\psi \: dx,\qquad \forall
\psi \in C_0^\infty (D_0),
\end{equation}
whence it follows that there exists a weak derivative $\partial_{x_2}u_0^{-}$ and $\gamma_2 = h\, \partial_{x_2}u_0^{-}$ a. e. in $D_0.$

Now let us find $\gamma _1.$ Consider the function
$$ \Phi(x) =
\begin{cases}
   0,& x\in \Omega _0, \\
 \varepsilon Y_1\left(\frac{x_1}{\varepsilon }\right)\psi + g, & x\in G_\varepsilon ,
  \end{cases}
\qquad \forall \psi \in C_0^\infty(D_0), \quad \psi\geq0,
$$
where  $Y_1(\xi  ) =-\xi +[\xi ]$. It is easy to see that $\Phi \in
K_\varepsilon $ and
$$
\nabla (\Phi -g) =\left(-\psi +\varepsilon
Y_1\left(\frac{x_1}{\varepsilon }\right)\partial _{x_1}\psi ,\,\,
\varepsilon Y_1\left(\frac{x_1}{\varepsilon }\right)\partial
_{x_2}\psi \right),\qquad x\in G_\varepsilon .
$$
Substituting the function $\Phi-g$ into the integral inequality
(\ref{2.4}) for solution $u_{\varepsilon},$ we get
\begin{equation*}
\int\limits_{G_\varepsilon }\left(- \partial_{x_1}u_\varepsilon\, \psi  + \varepsilon
Y_1\left(\frac{x_1}{\varepsilon }\right) \nabla u_\varepsilon\cdot \nabla \psi \right)\: dx\geq
\varepsilon \int\limits_{G_\varepsilon} Y_1\left(\frac{x_1}{\varepsilon }\right) f \,\psi
\:dx - \frac{\varepsilon^2(1\pm h)}{2}\int\limits_{S^{\pm}_\varepsilon }d \, \psi \:dx_2,
\end{equation*}
where the sings $"+"$ or $"-"$ in $S^{\pm}_\varepsilon$ indicate the union of the right or left sides of the thin rods respectively. With the help of (\ref{4.2}) and (\ref{4.7}) we deduce from previous inequality the estimate
\begin{equation*}
\begin{split}
 \left|\int_{D_0}\widetilde{\partial_{x_1} u_\varepsilon}\,\psi \:dx\right| & \leq \varepsilon \left(
\int_{G_\varepsilon }\bigl|Y_1\left(\frac{x_1}{\varepsilon
}\right)\left( \nabla u_\varepsilon \cdot\nabla \psi - f \,\psi
\right)\bigr|dx + \frac{\varepsilon (1+h)}{2} \int_{S_\varepsilon
}|d \, \psi |\: dx_2\right)\leq
\\[1 mm]
 \leq &\,  \varepsilon c_1
\left(\|\nabla u_\varepsilon \|_{L^2(G_\varepsilon)} \|\nabla \psi \|_{L^2(G_\varepsilon )} +
\|f\|_{L^2(G_\varepsilon )} \|\psi\|_{L^2(G_\varepsilon)} +
\varepsilon \| d \|_{L^2\left(S_\varepsilon\right)} \|\psi \|_{L^2\left(S_\varepsilon\right)}\right)\leq
\\[1 mm]
 \leq &  \, \varepsilon c_1\left(\|u_\varepsilon\|_{H^1(\Omega_\varepsilon)}
\|\psi \|_{H^1(D_0)} +  \|f\|_{L^2(\Omega_1)} \|\psi \|_{L^2(D_0)} + \|d\|_{H^1(D_0)}  \|\psi\|_{H^1(D_0)}\right)\leq
\varepsilon c_2,
\end{split}
\end{equation*}
from which, passing to the limit as $\varepsilon \rightarrow 0$, we get
$\int_{D_0}\gamma _1\psi \:dx=0$ for all $\psi \in C_0^\infty (D_0),\  \psi\geq0.$
This means that $\gamma_1=0$ a. e. in $D_0.$

\paragraph{2.} Let us show that the traces of the functions $u_0^+$ and $u_0^{-}$ on $I_0$ are equal.
By virtue of the compactness of the trace operator and the first
relation in (\ref{5.4}), we have
\begin{equation}\label{5.5}
u_\varepsilon(x_1,0) \stackrel{s}{\longrightarrow}u_0^+(x_1,0)\quad
\mbox{in}\ L^2(0,a)\quad\mbox{as}\quad \varepsilon \rightarrow 0.
\end{equation}

Consider the following equality
\begin{equation}\label{5.6}
\widetilde{u_\varepsilon}(x_1,0)=\chi
_{h}\left(\frac{x_1}{\varepsilon }\right)u_\varepsilon (x_1,0) \quad
\mbox{for a. e. }x_1\in (0,a),
\end{equation}
where $\chi_{h}(\xi ),\:\:\xi \in\mathbb{R}$, is the 1-periodic
function, defined on the  segment $[0,1]$ as follows:
$$
\chi _{h}(\xi )= \begin{cases}
    1, & |\xi -\frac{1}{2}|\leq\frac{h}{2}, \\
    0, & \frac{h}{2}<|\xi -\frac{1}{2}|\leq 1.
  \end{cases}
$$
It is known that $\chi _{h}\left(\frac{x_1}{\varepsilon
}\right)\stackrel{w}{\longrightarrow} h$ \ weakly in \ $L^2(0,1)$ \ as \ $ \varepsilon \rightarrow 0.$ Using this
fact and (\ref{5.5}), we obtain that the right-hand side in (\ref{5.6}) converges to $h \,
u_0^+(x_1,0)$ weakly in $L^2(0,a).$
 On the other hand,
\begin{equation}\label{5.7}
\int\limits_{0}^a\widetilde{u_\varepsilon} (x_1,0)\psi
(x_1)\:dx_1=\frac{1}{l}\int\limits_{D_0}\widetilde{u_\varepsilon}(x)
\psi (x_1)\:dx + \frac{1}{l}\int\limits_{D_0} (x_2+l)\widetilde{\partial
_{x_2}u_\varepsilon} \cdot\psi (x_1)\:dx \quad \forall\: \psi
\in C_0^\infty(0,a).
\end{equation}
Passing to the limit in (\ref{5.7}) as $\varepsilon \to 0$ and taking (\ref{gamma_2}) into account, we have
\begin{equation*}
h\int\limits_{0}^a u_{0}^{+} (\cdot,0)\psi
(x_1)\:dx_1=\frac{h}{l}\int\limits_{D_0} u_{0}^{-} \psi (x_1)\:dx +
\frac{h}{l}\int\limits_{D_0} (x_2+l)\, \partial _{x_2}u_{0}^{-}\, \psi (x_1)\:dx \quad \forall\: \psi \in C_0^\infty(0,a),
  \end{equation*}
whence it appears
$$\int_{0}^a u_{0}^{+}(\cdot,0)\psi (x_1)\:dx_1=\int_{0}^a u_{0}^{-} (\cdot,0)\psi
(x_1)\:dx_1 \quad \forall\: \psi \in C_0^\infty(0,a),
$$
i.e., $u_0^+(x_1,0)=u_0^{-}(x_1,0)$ for a.e. $x_1\in (0,a).$

Similarly we can prove that the trace $u_0^{-}|_{I_l}$ is equal to zero.

Thus, the results obtained above mean that the function
$$
u_0(x) =\left\{
\begin{array}{ll}
  u_0^+ , & x\in  \Omega_0 , \\
  u_0^-,  & x\in D_0 ,
\end{array}\right.
$$
belongs to the space ${\cal H}(\Omega_1, D_0; I_l).$

\paragraph{3.}
Let us add the inequality
$$
\int_{\Omega_0}\nabla (\varphi-u_{\varepsilon})\cdot \nabla
(\varphi-u_{\varepsilon})
\,dx+\int_{G_{\varepsilon}}\partial_{x_2}(\varphi-u_{\varepsilon})
\,\partial_{x_2}(\varphi-u_{\varepsilon})\, dx +\int_{G_{\varepsilon}}
\partial_{x_1}u_{\varepsilon} \, \partial_{x_1} u_{\varepsilon}\, dx\geq0,
$$
where $ \varphi$ is arbitrary function from $C^1(\overline{\Omega_1})$ such that
$\left.\varphi\right|_{I_l}=0$ and $\varphi\leq g$ in $D_0$ (obviously $\varphi|_{\Omega_{\varepsilon}}\in
K_\varepsilon ),$ \ to inequality (\ref{ozn2.5}). We get
$$
\int_{\Omega_0}\nabla \varphi \cdot \nabla (\varphi -
u_{\varepsilon})\,dx + \int_{G_{\varepsilon}}\partial_{x_1} u_{\varepsilon} \, \partial_{x_1} \varphi \,dx
+\int_{G_{\varepsilon}}\partial_{x_2} \varphi \, \partial_{x_2} (\varphi - u_{\varepsilon})\,dx\geq
$$
\begin{equation*}
\geq\int_{\Omega_{\varepsilon}} f (\varphi-u_{\varepsilon})
\,dx +\varepsilon \int_{ S_{\varepsilon}} d(x)
(\varphi-u_{\varepsilon}) \,ds,
\end{equation*}
which we can re-write with the help of (\ref{4.1}) in the following view
\begin{multline}\label{5.8}
 \int\limits_{\Omega_{0}}\nabla
\varphi\cdot\nabla(\varphi-u_{\varepsilon})\:dx
+\int\limits_{D_0}\widetilde{\partial_{x_1} u_{\varepsilon}} \,
\partial_{x_1} \varphi \,dx
+\int\limits_{D_0} \chi _{h}\left(\frac{x_1}{\varepsilon }\right)
\partial_{x_2} \varphi \,
\partial_{x_2} \varphi \, dx - \int\limits_{D_0}\partial_{x_2} \varphi \, \widetilde{\partial_{x_2}
 u_{\varepsilon}}\,dx\ \geq
\\
 \geq \int\limits_{\Omega_0}f \, (\varphi -u_{\varepsilon})\:
dx+\int\limits_{D_0}\chi _{h}\left(\frac{x_1}{\varepsilon }\right)f
\, \varphi \: dx  -\int\limits_{D_0}f \widetilde{u_{\varepsilon}} \,
dx +
\\
+
\frac{2}{h}\int\limits_{D_0}\chi _{h}\left(\frac{x_1}{\varepsilon }\right) d\,\varphi \: dx -
 \frac{2}{h}\int\limits_{D_0} d\, \widetilde{u_{\varepsilon}} \: dx
- \frac{2\varepsilon}{h}\int_{G_{\varepsilon} }Y\left(\frac{x_1}{\varepsilon
}\right)\partial_{x_1}(d (\varphi-u_{\varepsilon})) \: dx.
\end{multline}

Passing to the limit in (\ref{5.8}) as $\varepsilon
\rightarrow0$ and taking into account results obtained above, we obtain the following integral
inequality
$$
\int\limits_{\Omega_{0}}\nabla \varphi \cdot \nabla (\varphi-u_0^+)
\, dx+h \int\limits_{D_0} \partial_{x_2}\varphi\,
\partial_{x_2}(\varphi - u_0^-)\, dx\geq
$$
\begin{equation}\label{5.9}
\geq\int\limits_{\Omega_0}f (\varphi-u_0^+) \, dx +h
\int\limits_{D_0} f \,(\varphi -u_0^-)\, dx+2 \int\limits_{D_0} d
\,(\varphi -u_0^-)\, dx
\end{equation}
for any function
$\varphi \in K_1= \left\{\varphi \in C^1(\overline{\Omega_1}):
\ \varphi\left.\right|_{I_l}=0, \ \varphi\leq g \
\mbox{in} \ D_0\right\}.$

Since the set $K_1$ is dense in $K_0,$ the integral inequality (\ref{5.9}) holds for any function $\varphi\in K_0.$
This means that the function $u_0$ is the unique solution of inequality (\ref{def4.1}) (see Definition~\ref{def-4})
and also it is the weak solution to the homogenized problem (\ref{limit-problem}).

Due to the uniqueness of the solution to problem (\ref{limit-problem}), the above argumentations are true  for any subsequence of $\{\varepsilon\}$ chosen at the beginning of the proof. Thus the limits (\ref{limits}) hold.

\paragraph{4.}
>From equalities (\ref{2.1}) and (\ref{def3.1}) it follows that
\begin{equation}\label{5.10}
E_\varepsilon(u_\varepsilon)=\int\limits_{\Omega_{\varepsilon}}|\nabla
u_{\varepsilon}|^2\,dx= \int\limits_{G_{\varepsilon}}\nabla
u_{\varepsilon}\cdot \nabla g\,dx+ \int\limits_{\Omega_{0}} f
u_{\varepsilon}\,dx + \int\limits_{G_{\varepsilon}} f
(u_{\varepsilon}-g)dx+\varepsilon\int\limits_{ S_{\varepsilon}}
d(x)(u_{\varepsilon}-g)\, ds,
\end{equation}
$$
E_0(u_0)=\int\limits_{\Omega_{0}}|\nabla u_0^+|^2 \, dx+h
\int\limits_{D_0} |\partial_{x_2}u_0^-|^2\, dx=h\int\limits_{D_0}
\partial_{x_2}u_0^- \,\partial_{x_2}g\, dx+
$$
\begin{equation}\label{5.11}
 +\int\limits_{\Omega_0}f u_0^+ \, dx +h
\int\limits_{D_0} f \,(u_0^- -g)\, dx+2 \int\limits_{D_0} d \,(u_0^-
-g)\, dx .
\end{equation}
Passing to the limit in (\ref{5.10}) similarly as we made this in (\ref{5.8})
and taking into account  (\ref{5.11}), we obtain
\ $\lim_{\varepsilon \to 0} E_{\varepsilon}(u_{\varepsilon})= E_{0}(u_0).$ The theorem is proved.


\begin{thebibliography}{99}
\bibitem{B-Gaudiello-Melnyk}
{\sc D.~Blanchard, A.~Gaudiello, and T.~A.~Mel'nyk},
{\em Boundary homogenization and reduction of dimention in a
Kirchhoff-Love plate}, SIAM Journal on Mathematical
Analysis, {\bf 39} (2008), no. 6, pp. 1764--1787.

\bibitem{B-Gaudiello-Mos}
{\sc D.~Blanchard, A.~Gaudiello, and J.~Mossino},  {\em Highly
oscillating boundaries and  reduction of dimension: the critical
case}, Analysis and Application,  {\bf 5} (2007), pp.~137--163.

\bibitem{B-Gaudiello-Mos1}
{\sc D.~Blanchard, A.~Gaudiello, and G. Griso}, {\em Junction of a periodic family of elastic rods with 3d plate. Part I.}, J. Math. Pures Appl., {\bf 88} (2007), no. 9, 1-33 (Part I); {\bf 88} (2007), no. 9, 149-190 (Part II).

\bibitem{Chechkin-Melnyk-1}
{\sc G.A.~Chechkin, T.A.~Mel'nyk}, {\em Asymptotic analysis of boundary-value problems in thick
cascade junctions}, Reports of National Ukrainian Academy of Sciences, {\bf 9}, (2008),

\bibitem{D-D-M-07}
{\sc C.~D'Apice, U.~De Maio, and T.~A. Mel'nyk}, {\em Asymptotic
analysis of a perturbed parabolic problem in a thick junction of
type 3:2:2}, Networks and Heterogeneous Media, {\bf 2} (2007), 255--277.

\bibitem{Mel-m2as-08}
{\sc T.A. Mel'nyk}, {\em Homogenization of a boundary-value problem with
a nonlinear boundary condition in a thick junction of type 3:2:1},
Mathematical Models and Methods in Applied Sciences, {\bf 31} (2008), no.~9, 1005--1027.
Published online:\ {\tt http://dx.doi.org/10.1002/mma.951}

\bibitem{M-V-1}
{\sc T.A.~Mel'nyk, P.S.~Vaschuk},  {\em
Homogenization of the Neumann-Fourier Problem in a Thick Two-level Junction
of Type 3:2:1}, J. of Math. Physics, Analysis and Geometry,
{\bf 2} (2006), 318--337.

\bibitem{M-V-2}
{\sc T.A.~Mel'nyk,  P.S.~Vaschuk},  {\em
Homogenization of a boundary-value problem with mixed type of boundary conditions in a thick junction},
Differential Equations, {\bf 43} (2007), no. 5,  696--703.

\bibitem{UMJ}
{\sc T.A.~Mel'nyk}, {\em Homogenization of a singularly perturbed
parabolic problem in a thick periodic junction of the type 3:2:1},
Ukrainskii Matem. Zhurnal, \textbf{52} (2000), pp. 1524--1534 [in Ukrainian];
English transl.: Ukrainian Math. Journal, \textbf{52} (2000), 1737--1749.

\bibitem{M-steklov}
{\sc T.A. Mel'nyk}, {\em Asymptotic behavior of eigenvalues
and eigenfunctions of the Steklov problem in a thick periodic
junction}, Nonlinear oscillations, \textbf{4} (2001), no. 1,  91--105.

\bibitem{MN94}
{\sc T.A. Mel'nyk, S.A. Nazarov}, {\em Asymptotic structure of
the spectrum of the Neumann problem in a thin comb-like domain}, C.R. Acad
Sci. Paris, Serie 1, \textbf{319} (1994), 1343--1348.

\bibitem{N95}
{\sc S.A. Nazarov}, {\em Junctions of singularly degenerating
domains with different limit dimensions}, Trudy  Seminara imeni I.G. Petrovskogo,
\textbf{18} (1995), 3-79 (part I);  \textbf{20} (2000),
155--196 (part II).

\bibitem{MN96}
{\sc T.A. Mel'nyk,  S.A. Nazarov}, {\em Asymptotics of the
Neumann spectral problem solution in a domain of ``thick comb"
type}, Trudy Seminara imeni I.G. Petrovskogo, \textbf{19} (1996)
138--173 [in Russian]; English transl.: J. of Math. Sciences, \textbf{85} (1997), 2326--2346.

\bibitem{ZAA99}
{\sc T.A. Mel'nyk}, {\em Homogenization of the Poisson equation
in a thick periodic junction}, Zeitschrift f\"ur Analysis und ihre
Anwendungen, \textbf{18} (1999),  953--975.

\bibitem{M00}
{\sc T.A. Mel'nyk}, {\em Asymptotic analysis of a spectral problem
in a periodic thick junction of type 3:2:1}, Mathematical Methods in
the Applied sciences, \textbf{23} (2000),  321--346.

\bibitem{M-N-AiA}
{\sc T.A. Mel'nyk,  S.A. Nazarov}, {\em Asymptotic analysis
of the Neumann problem of the junction of a body and thin heavy
rods}, Algebra i Analiz, \textbf{12} (2000),  188--238 [in Russian];
English transl.: St.Petersburg Math.J., \textbf{12}, (2001),
317--351.

\bibitem{Fleury-SanchezPal}
{\sc F. Fleury,  E. Sanchez-Palencia}, {\em Asymptotic and
spectral properties of the acoustic vibrations of body, perforated by narrow
channels}, Bull. Sci. Math., \textbf{2} (1986), 149--176.

\bibitem{Khrus}
{\sc E.Ya. Khruslov}, {\em On the resonance phenomenas in the one
problem of diffraction}, Teorija Funktsij, Funcionalnyj Analis i
ix Prilozhenija (Izd-vo Kharkov Univ.), {\bf 8} (1968),  113--120
(in Russian).

\bibitem{Khrus-Kot}
{\sc V.P. Kotliarov, E.Ya. Khruslov}, {\em On a limit boundary condition
of some Neumann problem},  Theor. Funkts., Funkts. Anal.
Prilozhen (Izd-vo Kharkov Univ.), {\bf 10} (1970), 83-96 (in Russian).

\bibitem{Signorini}
{\sc A. Signorini}, {\em Questioni di elasticita non linearizzata o
semilinearizzata}, Rend. di Matem. e delle sue appl., {\bf 18} (1959).

\bibitem{Kinder-Stam}
{\sc D. Kinderlehrer, G. Stampaccia},  {\em An introduction to variational inequalities and their applications}.
Academic Press, 1980.

\bibitem{Lions-Stam}
{\sc J.-L. Lions, G. Stampaccia}, {\em Variational inequalities}, Comm. Pure Appl. Math. {\bf 20} (1976), 493--519.

\bibitem{Lions}
{\sc J.-L. Lions}, {\em Quelques m\'{e}thodes de r\'{e}solution des probl$\rm\grave{e}$mes aux limites non lin\'{e}ires}.  Dunod, Paris, 1969.

\bibitem{Zhikov-1}
{\sc V.V. Zhikov}, {\em On the homogenization of nonlinear variational problems in perforated domains}, Russian Journal of Mathematical Physics. {\bf 12} (1994) no.~3, 393--408.

\bibitem{Pastukhova}
{\sc S. E. Pastukhova}, {\em Homogenization of a mixed problem with
Signorini condition for an  elliptic operator  in a perforated
domain}, Sb. Math. {\bf 192} (2001), 245-260.

\bibitem{Vor-Shaposhnikova}
{\sc A.Yu. Vorobev, T.S. Shaposhnikova}, {\em On homogenization of the nonuniform Signorini problem for the Poisson equation in a periodic perforated domain}, Differential Equations, {\bf 39} (2003), no. 3, 359-366 (in Russian).

\bibitem{Sandrakov-1}
{\sc G.V. Sandrakov}, {\em Homogenization of variational inequalities for
problems with regular obstacles}, Dokl. Akad. Nauk, {\bf 397} (2004),
170-173; English transl., Dokl. Math. {\bf 71} (2004), 119-122.

\bibitem{Sandrakov-2}
{\sc G.V. Sandrakov}, {\em Homogenization of variational inequalities for non-linear diffusion problems in perforated domains}. Izvestia: Mathematics, {\bf 69} (2005), no. 5, 1035-1059.

\bibitem{Mel-Kaz}
{\sc Yu.A. Kazmerchuk, T.A. Mel'nyk}, {\em Homogenization of the Signorini boundary-value\\ problem
in a thick plane junction}, Nonlinear Oscillations, {\bf 12} (2009), no. 1 (to appear).

\bibitem{U}
{\sc S.V. Uspenskii}, {\em The traces of functions the Sobolev
space $W^{l_1,\ldots,l_n}_p$ on smooth surfaces}. Siberian
Math. J. {\bf 13} (1972), 298--313.

\end{thebibliography}
\end{document}